\documentclass[12pt,amsmath]{article}
\usepackage{amsfonts}
\begin{document}
\newsymbol\rtimes 226F
\newsymbol\ltimes 226E
\newcommand{\text}[1]{\mbox{{\rm #1}}}
\newcommand{\Rep}{\text{Rep}}
\newcommand{\gr}{\text{gr}}
\newcommand{\Fun}{\text{Fun}}
\newcommand{\Hom}{\text{Hom}}
\newcommand{\End}{\text{End}}
\newcommand{\FPdim}{\text{FPdim}}
\newcommand{\GL}{\text{GL}}
\newcommand{\Sp}{\text{Sp}}
\newcommand{\Ps}{\text{Ps}}
\newcommand{\Ad}{\text{Ad}}
\newcommand{\ASp}{\text{ASp}}
\newcommand{\APs}{\text{APs}}
\newcommand{\Rad}{\text{Rad}}
\newcommand{\Corad}{\text{Corad}}
\newcommand{\SuperVect}{\text{SuperVect}}
\newcommand{\Vect}{\text{Vect}}
\newcommand{\Spec}{\text{Spec}}
\newcommand{\tr}{\text{tr}}
\newcommand{\cH}{{\cal H}}
\newcommand{\cR}{{\cal R}}
\newcommand{\cJ}{\cal J}
\newcommand{\gd}{\delta}
\newcommand{\lan}{\langle}
\newcommand{\ran}{\rangle}
\newcommand{\itms}[1]{\item[[#1]]}
\newcommand{\nin}{\in\!\!\!\!\!/}
\newcommand{\g}{{\bf g}}
\newcommand{\sub}{\subset}
\newcommand{\cntd}{\subseteq}
\newcommand{\go}{\omega}
\newcommand{\Pa}{P_{a^\nu,1}(U)}
\newcommand{\fx}{f(x)}
\newcommand{\fy}{f(y)}
\newcommand{\gD}{\Delta}
\newcommand{\gl}{\lambda}
\newcommand{\gL}{\Lambda}
\newcommand{\half}{\frac{1}{2}}
\newcommand{\sto}[1]{#1^{(1)}}
\newcommand{\stt}[1]{#1^{(2)}}
\newcommand{\Z}{\hbox{\sf Z\kern-0.720em\hbox{ Z}}}
\newcommand{\singcolb}[2]{\left(\begin{array}{c}#1\\#2
\end{array}\right)}
\newcommand{\ga}{\alpha}
\newcommand{\gb}{\beta}
\newcommand{\gga}{\gamma}
\newcommand{\ul}{\underline}
\newcommand{\ol}{\overline}
\newcommand{\qed}{\kern 5pt\vrule height8pt width6.5pt depth2pt}
\newcommand{\Lrraro}{\Longrightarrow}
\newcommand{\Nb}{|\!\!/}
\newcommand{\NN}{{\rm I\!N}}
\newcommand{\bsl}{\backslash}
\newcommand{\gt}{\theta}
\newcommand{\op}{\oplus}
\newcommand{\C}{{\bf C}}
\newcommand{\Q}{{\bf Q}}
\newcommand{\Op}{\bigoplus}
\newcommand{\CR}{{\cal R}}
\newcommand{\grr}{\omega_1}
\newcommand{\ben}{\begin{enumerate}}
\newcommand{\een}{\end{enumerate}}
\newcommand{\ndiv}{\not\mid}
\newcommand{\bab}{\bowtie}
\newcommand{\hal}{\leftharpoonup}
\newcommand{\har}{\rightharpoonup}
\newcommand{\ot}{\otimes}
\newcommand{\OT}{\bigotimes}
\newcommand{\bwe}{\bigwedge}
\newcommand{\eps}{\varepsilon}
\newcommand{\gs}{\sigma}
\newcommand{\rbraces}[1]{\left( #1 \right)}
\newcommand{\bbox}{$\;\;\rule{2mm}{2mm}$}
\newcommand{\sbraces}[1]{\left[ #1 \right]}
\newcommand{\bbraces}[1]{\left\{ #1 \right\}}
\newcommand{\OO}{_{(1)}}
\newcommand{\TT}{_{(2)}}
\newcommand{\FF}{_{(3)}}
\newcommand{\minus}{^{-1}}
\newcommand{\CV}{\cal V}
\newcommand{\CVs}{\cal{V}_s}
\newcommand{\un}{U_q(sl_n)'}
\newcommand{\on}{O_q(SL_n)'}
\newcommand{\slq}{U_q(sl_2)}
\newcommand{\olq}{O_q(SL_2)}
\newcommand{\UU}{U_{(N,\nu,\go)}}
\newcommand{\HH}{{\mathcal H}}
\newcommand{\ct}{\centerline}
\newcommand{\bs}{\bigskip}
\newcommand{\qua}{\rm quasitriangular}
\newcommand{\ms}{\medskip}
\newcommand{\noin}{\noindent}
\newcommand{\mat}[1]{$\;{#1}\;$}
\newcommand{\raro}{\rightarrow}
\newcommand{\map}[3]{{#1}\::\:{#2}\raro{#3}}
\newcommand{\alg}{{\rm Alg}}
\def\newtheorems{\newtheorem{theorem}{Theorem}[subsection]
                 \newtheorem{cor}[theorem]{Corollary}
                 \newtheorem{proposition}[theorem]{Proposition}
                 \newtheorem{lemma}[theorem]{Lemma}
                 \newtheorem{defn}[theorem]{Definition}
                 \newtheorem{Theorem}{Theorem}[section]
                 \newtheorem{Corollary}[Theorem]{Corollary}
                 \newtheorem{corollary}[theorem]{Corollary}
                 \newtheorem{Proposition}[Theorem]{Proposition}
                 \newtheorem{Lemma}[Theorem]{Lemma}
                 \newtheorem{Definition}[Theorem]{Definition}
                 \newtheorem{Example}[Theorem]{Example}
                 \newtheorem{Remark}[Theorem]{Remark}
                 \newtheorem{claim}[theorem]{Claim}
                 \newtheorem{sublemma}[theorem]{Sublemma}
                 \newtheorem{example}[theorem]{Example}
                 \newtheorem{definition}[theorem]{Definition}
                 \newtheorem{remark}[theorem]{Remark}
                 \newtheorem{question}[theorem]{Question}
                 \newtheorem{Question}[Theorem]{Question}
                 \newtheorem{Conjecture}[Theorem]{Conjecture}}
\newtheorems
\newcommand{\proof}{\par\noindent{\bf Proof:}\quad}
\newcommand{\dmatr}[2]{\left(\begin{array}{c}{#1}\\
                            {#2}\end{array}\right)}
\newcommand{\doubcolb}[4]{\left(\begin{array}{cc}#1&#2\\
#3&#4\end{array}\right)}
\newcommand{\qmatrl}[4]{\left(\begin{array}{ll}{#1}&{#2}\\
                            {#3}&{#4}\end{array}\right)}
\newcommand{\qmatrc}[4]{\left(\begin{array}{cc}{#1}&{#2}\\
                            {#3}&{#4}\end{array}\right)}
\newcommand{\qmatrr}[4]{\left(\begin{array}{rr}{#1}&{#2}\\
                            {#3}&{#4}\end{array}\right)}
\newcommand{\smatr}[2]{\left(\begin{array}{c}{#1}\\
                            \vdots\\{#2}\end{array}\right)}

\newcommand{\ddet}[2]{\left[\begin{array}{c}{#1}\\
                           {#2}\end{array}\right]}
\newcommand{\qdetl}[4]{\left[\begin{array}{ll}{#1}&{#2}\\
                           {#3}&{#4}\end{array}\right]}
\newcommand{\qdetc}[4]{\left[\begin{array}{cc}{#1}&{#2}\\
                           {#3}&{#4}\end{array}\right]}
\newcommand{\qdetr}[4]{\left[\begin{array}{rr}{#1}&{#2}\\
                           {#3}&{#4}\end{array}\right]}

\newcommand{\qbracl}[4]{\left\{\begin{array}{ll}{#1}&{#2}\\
                           {#3}&{#4}\end{array}\right.}
\newcommand{\qbracr}[4]{\left.\begin{array}{ll}{#1}&{#2}\\
                           {#3}&{#4}\end{array}\right\}}

\title{{\bf Classification of Finite-Dimensional Triangular
Hopf Algebras With The Chevalley Property}}
\author{
Pavel Etingof $^1$ \and Shlomo Gelaki $^2$ } \footnotetext[1]
{MIT, Department of Mathematics, 77 Massachusetts Avenue,
Cambridge, MA 02139, USA, \& Columbia University, Department of
Mathematics, 2990 Broadway, New York, NY 10027, USA, {\rm email:
etingof@math.mit.edu} } \footnotetext[2] {Technion-Israel
Institute of Technology, Department of Mathematics, Haifa 32000,
Israel, {\rm email:}\linebreak {\rm gelaki@math.technion.ac.il} }

\maketitle

\section{Introduction}

Recall [AEG] that a triangular Hopf algebra $A$ over $\C$ is said
to have the Chevalley property if the tensor product of any two
simple $A$-modules is semisimple, or, equivalently, if the radical
of $A$ is a Hopf ideal. There are two reasons to study this class
of triangular Hopf algebras: First, it contains all known examples
of finite-dimensional triangular Hopf algebras; second, it can be,
in a sense, "completely understood". Namely, it was shown in [AEG]
that any finite-dimensional triangular Hopf algebra with the
Chevalley property is obtained by twisting of a finite-dimensional
triangular Hopf algebra with $R-$matrix of rank $\le 2,$ which, in
turn, is obtained by "modifying" the group algebra of a finite
supergroup. This provides a classification of such Hopf algebras.

The goal of this paper is to make this classification more
effective and explicit, i.e. to parameterize isomorphism classes
of finite-dimensional triangular Hopf algebras with the Chevalley
property by group-theoretical objects, similarly to how it was
done in [EG2] in the semisimple case. This is achieved in Theorem
\ref{main1}, where these classes are put in bijection with certain
septuples of data. In the semisimple case, the septuples reduce to
the quadruples of [EG2], and we recover the result of [EG2]. In
the minimal triangular pointed case, we recover Theorem 5.1 of
[G].

{\bf Acknowledgments.} The second author is grateful to MIT for
its hospitality. The authors were partially supported by the NSF
grant DMS-9988796. The work of the first author was done in part
for the Clay Mathematics Institute.

\section{The main theorem}

In this section we give an explicit description of the set of
isomorphism classes of finite-dimensional triangular Hopf algebras
over $\C$ with the Chevalley property. We will freely use the
facts from the theory of Hopf superalgebras and finite supergroups
which were sketched in [AEG].

\begin{Definition}\label{ept}
A triangular septuple is a
septuple $(G,W,H,Y,B,V,u)$ where $G$ is a finite group, $W$
is a finite-dimensional representation of $G,$ $H$ is a subgroup of $G,$
$Y$ is an $H-$invariant subspace of $W,$ $B$ is an
$H-$invariant nondegenerate element in $S^2Y,$ $V$ is an irreducible
projective representation of $H$ of dimension $|H|^{1/2},$ and $u\in G$
is a central element of order $\le 2$ acting by $-1$ on $W$.
\end{Definition}

The notion of isomorphism of triangular septuples is obvious.

Given a triangular septuple, one can construct a
finite-dimensional triangular Hopf algebra $A(G,W,H,Y,B,V,u)$ as
follows.

Regard $Y$ as a purely odd supervector space and consider the
supergroup $H\ltimes Y.$ Consider the group algebra $\C[H\ltimes
Y]$ of this supergroup. Let $J_V$ be a (minimal) twist for $\C[H]$
corresponding to $(H,V)$ as in [EG2] (it is well defined only up
to gauge equivalence). Let $J_B:=e^{B/2},$ and define ${\cal
J}:=J_BJ_V$. Then ${\cal J}$ is a twist for $\C[H\ltimes Y],$ and
it is clear (since $B$ is invariant) that the gauge equivalence
class of ${\cal J}$ is independent of the choice of $J_V$.

Regard $W$ as an odd supervector space, and consider the supergroup
$G\ltimes W$. We have a natural inclusion of supergroups $H\ltimes
Y\hookrightarrow G\ltimes W,$ hence ${\cal J}$ can be regarded as
a twist for the supergroup algebra $\C[G\ltimes W]$. Now let
$A(G,W,H,Y,B,V,u)$ be the finite-dimensional triangular Hopf
algebra corresponding to the pair $(\C[G\ltimes W]^{\cal J},u)$
under the correspondence of Theorem 3.3.1 of [AEG]. According to
[AEG], this Hopf algebra has the Chevalley property.

Our main result is the following theorem.

\begin{Theorem}\label{main1}

The assignment $(G,W,H,Y,B,V,u)\to A(G,W,H,Y,B,V,u)$
is a bijection between:

\ben

\item isomorphism classes of triangular septuples, and

\item isomorphism classes of finite-dimensional triangular Hopf algebras
over $\C$ with the Chevalley property.

\een

\end{Theorem}

\begin{Remark} {\rm It is clear that $A(G,W,H,Y,B,V,u)$ is semisimple
if and only if $W=0$ (and hence $Y=0$ and $B=0$). Therefore,
Theorem \ref{main1} implies, in particular,
 that there is a natural bijection between
isomorphism classes of semisimple triangular Hopf algebras over
$\C,$ and isomorphism classes of quadruples $(G,H,V,u)$. This was
the main result of [EG2]. Thus, Theorem \ref{main1} is a
generalization of the main result of [EG2]. }
\end{Remark}

The rest of the section is devoted to the proof of Theorem \ref{main1}.

 We start with the
following proposition, which is essentially contained in [AEG].

\begin{Proposition}\label{p1} Any minimal triangular Hopf superalgebra
${\cal A}$ over $\C$ with Drinfeld element $1$ and the Chevalley
property is isomorphic to $\C[H\ltimes Y]^{\cal J},$ where $H$ is
a finite group, $Y$ is a finite-dimensional representation of $H$
(considered as an odd vector space) and ${\cal J}=e^{B/2}J,$ where
$B\in (S^2Y)^H,$ and $J$ is a minimal twist for $H.$
\end{Proposition}

\proof Let ${\cal A}_s:={\cal A}/\Rad({\cal A})$ be the semisimple
quotient of ${\cal A}.$ As was shown in [G] and [AEG], this is a
minimal triangular semisimple Hopf algebra with Drinfeld element
$1,$ and we have a sequence of Hopf superalgebra homomorphisms
${\cal A}_s\to {\cal A}\to {\cal A}_s$, with the composition being
the identity. Thus, by [EG1] and [EG2], ${\cal A}_s=\C[H]^{J}$
with $R-$matrix $J_{21}^{-1}J,$ where $H$ is a finite group,
 and $J$ is a minimal twist for it, corresponding to an
irreducible projective representation $V$ of $H$ of dimension
$|H|^{1/2}.$

Let us now consider the Hopf superalgebra ${\cal A}^{J^{-1}}.$ We
have a sequence of Hopf superalgebra homomorphisms
$\C[H]\rightarrow {\cal A}^{J^{-1}}\to \C[H],$ and the composition
is the identity. The projection of the $R$-matrix ${\cal R}$ of
${\cal A}^{J^{-1}}$ to $\C[H]\ot \C[H]$ is $1\ot 1,$ so it is
unipotent in $X\otimes Z$ for any two ${\cal A}-$modules $X,Z.$
Let $({\cal A}^{J^{-1}})_m$ be the minimal part of ${\cal
A}^{J^{-1}}.$ Since any $({\cal A}^{J^{-1}})_m-$ module is
contained in an ${\cal A}-$module (see e.g. [AEG]), we get that
${\cal R}$ is unipotent in any $X\ot Z,$ where $X,Z$ are $({\cal
A}^{J^{-1}})_m-$modules. But this means that $({\cal
A}^{J^{-1}})_m$ is local (since $({\cal A}^{J^{-1}})_m/\Rad(({\cal
A}^{J^{-1}})_m)$ is minimal triangular with $R$-matrix $1\ot 1,$
hence $1-$dimensional). Hence by [AEG], $({\cal A}^{J^{-1}})_m$ is
$\Lambda Y,$ where $Y$ is a vector space, with some triangular
structure $R_B:=e^B,$ where $B\in S^2Y$ is a nondegenerate
element. Thus taking $J':=e^{B/2},$ we get that $({\cal
A}^{J^{-1}})^{J'^{-1}}$ is supercocommutative, with triangular
structure $1\otimes 1$. Then Kostant's theorem [K] implies that
$({\cal A}^{J^{-1}})^{J'^{-1}}=\C[H\ltimes Y]$. Thus, ${\cal
A}=\C[H\ltimes Y]^{J'J}$ with triangular structure
$(J'J)_{21}^{-1}J'J.$ This concludes the proof of the proposition.
\qed

The following proposition will imply the converse of Proposition
\ref{p1}: For any $H,Y,B,J$ as above, the Hopf superalgebra ${\cal
A}:=\C[H\ltimes Y]^{e^{B/2}J}$ is minimal triangular.

\begin{Proposition}\label{p2}
If $({\cal A},{\cal R}_{{\cal A}})$ is a minimal triangular Hopf
superalgebra, $H$ is a group acting on ${\cal A}$ by Hopf
superalgebra automorphisms and $J$ is a minimal twist of $H,$ then
the Hopf superalgebra $((\C[H]\ltimes {\cal A})^J,J_{21}^{-1}{\cal
R}_{{\cal A}}J)$ is minimal triangular.
\end{Proposition}

\proof First note that since $\C[H]$ is a Hopf subalgebra of
$\C[H]\ltimes {\cal A},$ $J$ is a twist for $\C[H]\ltimes {\cal A},$
and hence ${\cal R}:=J_{21}^{-1}{\cal R}_{{\cal A}}J$ is a triangular
structure on $\C[H]\ltimes {\cal A}.$ So all we have to prove is that
${\cal R}$ is minimal.

Let $Z:=(\C[H]\ltimes {\cal A})_m$ be the minimal part, and
$f:=Id\ot \varepsilon: \C[H]\ltimes {\cal A}\to \C[H].$ It is
straightforward to check that $(Id\ot f)({\cal R})=J_{21}^{-1}J;$
the minimal $R-$matrix for $\C[H]^J.$ Hence $\C[H]\subset Z.$ This
implies that ${\cal R}_{{\cal A}}=J_{21}{\cal R}J^{-1}\in Z\ot Z$
(as $J\in Z\otimes Z$). So we get that ${\cal A} \subset Z$ as
well. Since $Z$ is an algebra, we conclude that it is equal to
$\C[H]\ltimes {\cal A}$ as desired. \qed

Now let us turn to the proof of Theorem \ref{main1}. Our job is to
prove the following two statements:

\begin{Proposition}\label{p4}
Any finite-dimensional triangular Hopf algebra $(A,R)$ over $\C$
with the Chevalley property is of the form $A({\cal S})$ for a
suitable triangular septuple ${\cal S}$.
\end{Proposition}

\begin{Proposition}\label{p5}
Let $A({\cal S}_1),A({\cal S}_2)$ be two isomorphic triangular
Hopf algebras. Then the septuples ${\cal S}_1$, ${\cal S}_2$ are
isomorphic.
\end{Proposition}

\noin {\bf Proof of Proposition \ref{p4}:} Let $(A,R)$ be a
finite-dimensional triangular Hopf algebra with the Chevalley
property. By [AEG], the Drinfeld element $u$ of $A$ satisfies
$u^2=1.$ Let $({\cal A},{\cal R})$ be the triangular Hopf
superalgebra corresponding to $(A,R)$ under the correspondence
given in [AEG, Theorem 3.3.1]. Recall that the Drinfeld element of
${\cal A}$ is equal to $1,$ and that $u$ acts on ${\cal A}$ by
parity. Clearly ${\cal A}$ has the Chevalley property as well.

Let ${\cal A}_m$ be the minimal part of ${\cal A}$. By Proposition
\ref{p1}, ${\cal A}_m=\C[H\ltimes Y]^{{\cal J}},$ where ${\cal
J}=e^{B/2} J.$ Hence ${\cal A}^{{\cal J}^{-1}}$ is cocommutative,
and by Kostant's theorem [K], ${\cal A}^{{\cal
J}^{-1}}=\C[G\ltimes W],$ where $G$ is a finite group containing
$H$ as a subgroup and $Y\subseteq W$ is an $H-$invariant
subrepresentation. Thus we have associated to $(A,R)$ a triangular
septuple ${\cal S}=(G,W,H,Y,B,V,u)$. It is clear that
$(A,R)=A(G,W,H,Y,B,V,u).$ The proposition is proved. \qed

\noin {\bf Proof of Proposition \ref{p5}:} Let ${\cal
S}_i=(G_i,W_i,H_i,Y_i,B_i,V_i,u_i)$, $i=1,2$, be two triangular
septuples, which yield isomorphic triangular Hopf algebras
$A_i:=A({\cal S}_i)$. We want to show that ${\cal S}_1,{\cal S}_2$
are isomorphic.

Let $f:A_1\to A_2$ be an isomorphism of triangular Hopf
algebras. The Drinfeld element of $A_i$ is $u_i$, so we have $f(u_1)=u_2$.

Let ${\cal A}_i$ be the Hopf superalgebra with Drinfeld element
$1$ corresponding to $A_i$ as in [AEG, Theorem 3.3.1]. Since
$f(u_1)=u_2$, we find that $f$ defines an isomorphism $f:{\cal
A}_1\to {\cal A}_2$ of triangular Hopf superalgebras. Thus, $f$
defines an isomorphism of their minimal parts: $f(({\cal
A}_1)_m)=({\cal A}_2)_m$. Hence, $f(\text{Corad}(({\cal
A}_1)_m))=\text{Corad}(({\cal A}_2)_m)$. But by Proposition
\ref{p2}, $({\cal A}_i)_m=\C[H_i\ltimes Y_i]^{e^{B_i/2}J_i}$, so
$\text{Corad}(({\cal A}_i)_m)=\C[H_i]^{J_i}$. We conclude that $f$
restricts to an isomorphism of triangular Hopf algebras $f:
\C[H_1]^{J_1}\to \C[H_2]^{J_2}$. Hence, $J_2f^{\otimes 2}
(J_1)^{-1}$ is a symmetric twist. So by [EG2], it is equal to
$\Delta(x)(x^{-1}\otimes x^{-1})$ for some invertible $x\in
\C[H_2],$ and we have $f=\text{Ad}(x)\circ \gamma$, where $\gamma$
is a group isomorphism $H_1\to H_2$.

Without loss of generality, we can assume that $x=1$ (otherwise we
can apply a gauge transformation by $x$ to switch to a situation
where $x=1$). In this case, we get $f^{\otimes 2}(J_1)=J_2$, which
implies $f^*(V_1)=V_2,$ where the left hand side is the usual
pullback.

Now, by definition, we have $f^{\otimes
2}((J_1)_{21}^{-1}e^{B_1}J_1)=(J_2)_{21}^{-1}e^{B_2}J_2$. Thus,
$f^{\otimes 2}(e^{B_1})=e^{B_2}$, i.e. $f^{\otimes 2}(B_1)=B_2$.
This means that $f(Y_1)=Y_2.$ This map has to be consistent with
the actions of $H_i$ on $Y_i$ and the map $f:H_1\to H_2$, since
$f$ is a homomorphism of algebras.

Finally, setting ${\cal J}_i:=e^{B_i/2}J_i$, we have $f^{\ot 2
}({\cal J}_1)= {\cal J}_2$. Therefore, $f$ defines an isomorphism
of triangular Hopf superalgebras $\C[G_1\ltimes W_1]\to
\C[G_2\ltimes W_2]$.

Summarizing, we see that $f$ sets up an isomorphism
between ${\cal S}_1$ and ${\cal S}_2$, as desired.
The proposition and the theorem are proved.
\qed

\section{Minimal triangular pointed Hopf algebras}

Here we apply Theorem \ref{main1} to classifying
minimal triangular Hopf algebras with the Chevalley property,
and minimal triangular pointed Hopf algebras.

\begin{Proposition}\label{p6}
The following hold:
\ben
\item The triangular Hopf algebra
$A=A(G,W,H,Y,B,V,u)$ is minimal if and only if $Y=W$ and $G$ is
generated by $H,u$. Thus, minimal triangular Hopf algebras over
$\C$ with the Chevalley property correspond to triangular
septuples ${\cal S}$ such that $W=Y$ and $G=<H,u>$.

\item The triangular Hopf algebra $A=A(G,W,H,Y,B,V,u)$ is minimal
pointed if and only if $Y=W$, $G$ is generated by $H$ and $u,$ and
$G$ is abelian. Thus, minimal triangular pointed Hopf algebras
over $\C$ correspond to triangular septuples ${\cal S}$ such that
$W=Y$, $G=<H,u>$, and $G$ is abelian. \een
\end{Proposition}

\proof We start by proving part 1. The ``only if'' direction is
clear. To prove the ``if'' direction, consider the minimal part of
the Hopf superalgebra ${\cal A}$ with Drinfeld element $1$,
corresponding to $A$. By Proposition \ref{p2}, it is $\C[H\ltimes
W]$ (as an algebra). The minimal part of $A$ is clearly generated
by $u$ and the minimal part of ${\cal A}$, so it is everything.

Part 2 follows easily from part 1, if we remember from [AEG] that
any minimal triangular pointed Hopf algebra has the Chevalley
property. \qed

\begin{Remark} {\rm Proposition \ref{p6} implies that for minimal
triangular Hopf algebras with the Chevalley property (unlike the
nonminimal ones), there are finitely many isomorphism classes in
any given dimension. Indeed, in a ``minimal'' triangular septuple,
the only ``continuous'' parameter is the tensor $B$, or,
equivalently, the $G$-invariant inner product $B^{-1}$ on $Y=W$.
But such a product, if exists, is clearly unique up to an
isomorphism. In particular, this implies that any minimal
triangular Hopf algebra with the Chevalley property is rigid, i.e.
does not have nontrivial deformations as a triangular Hopf
algebra. Since it is suspected that all finite-dimensional
triangular Hopf algebras over $\C$ have the Chevalley property
(see [AEG], Question 5.5.1), it would be interesting to check
whether any minimal triangular Hopf algebra is rigid. }
\end{Remark}

To conclude the paper, we recall that a classification of minimal
triangular pointed Hopf algebras was in fact obtained earlier in
[G] (the additional condition to be generated by grouplike and
skewprimitive elements imposed in [G] is in fact always satisfied,
as was shown in [AEG]). There, the isomorphism classes of such
Hopf algebras were parameterized by somewhat different
group-theoretical data than here. So, let us identify these two
kinds of data with each other. This, in particular, will give
another proof of Theorem 5.1 in [G], based on Theorem \ref{main1}.

Recall the definitions of the two types of the parameterizing
data.

{\bf Type 1 (this paper):} A $5-$tuple $(G,W,H,V,u)$, where $G$ is
a finite abelian group, $H$ a subgroup of $G$, $u$ an element of
$G$ of order $2$, $W$ a finite-dimensional 
representation of $G$ with an invariant
inner product, and $V$ a projective irreducible representation of
$H$ of dimension $|H|^{1/2}$, such that $G=<H,u>$, and $u|_W=-1$.

{\bf Type 2 ([G, Section 4]):} A triple $(G,\phi,n)$, where $G$ is
a finite abelian group, $\phi$ a nondegenerate skew-symmetric
bilinear form on $G$ with values in $\C^*$ (i.e.
$\phi(g,h)=\phi(h,g)^{-1}$), and $n: g\mapsto n_g$ a nonnegative
integer function defined on the set $I_\phi$ of elements $g\in G$
such that $\phi(g,g)=-1$, satisfying $n_g=n_{g^{-1}}$ (the other
data in [G], up to isomorphism, can be uniquely expressed via
$(G,\phi,n)$).

Below we will describe how to pass from data of type 1 to data of
type 2, and back.

{\bf 1. From Type 1 to Type 2:} The group $G$ is already given, so
we need to construct $\phi$ and $n.$ To construct $\phi$, recall
that the quadruple $(G,H,V,u)$ defines a minimal triangular
structure on $\C[G]$, which is nothing but a nondegenerate
skew-symmetric form on $G$; so let it be $\phi$. Finally, to
construct $n$, recall that the group $G$ is abelian, and is
identified with its dual via $\phi$, so the representation $W$ is
representable as a direct sum $\oplus_g n_g \phi(g,*)$, and we have
$n_{g^{-1}}=n_g$ since the representation is orthogonal.
Furthermore, since $\phi(u,g)=\phi(g,g)$, and $u|_W=-1$, we get
$\phi(g,g)=-1$ if $n_g>0$.

{\bf 2. From Type 2 to Type 1:} Again, the group $G$ is given, and
we need to construct the quadruple $(W,H,V,u)$. To construct
$(H,V,u)$, we take the quadruple $(G,H,V,u)$ which corresponds to
the triangular structure on $\C[G]$ defined by $\phi$. Finally, to
construct $W$, we set $W:=\oplus_g n_g\phi(g,*)$. It is clear that
this representation is orthogonal since $n_g=n_{g^{-1}}$, and that
$u$ acts by $-1$ in $W$ since $\phi(g,g)=-1$ whenever $n_g>0$.

It is clear that the constructed two correspondences are inverses
to each other. We leave it to the reader to show that the data of
type 1 and type 2 that correspond to each other under this rule,
define the same minimal triangular pointed Hopf algebra.

\end{document}